\documentclass[11pt]{amsart}
\usepackage{amssymb,latexsym,amsmath}
\input{xypic}
\xyoption{all}

\begin{document}

\begin{abstract} 
This article grew out of my talk in `The Legacy of Srinivasa Ramanujan' conference where I spoke about some  techniques to prove algebraicity results for the special values of symmetric cube $L$-functions attached to the Ramanujan $\Delta$-function. If one wishes to compare these different techniques, then one needs to compare various automorphic periods attached to the symmetric cube transfer of $\Delta.$ Motivated by this problem, in this article we provide comparison results for 
{\it Whittaker-Betti periods, Shalika-Betti periods} and {\it relative periods} attached to a given cohomological cuspidal automorphic representation of ${\rm GL}_{2n}$ over a totally real number field. 
\end{abstract}

\title[Automorphic Periods]{Comparison results for certain periods of \\
cusp forms on ${\rm GL}_{2n}$ over a totally real number field}

\author{\bf A. Raghuram}
\address{Indian Institute of Science Education and Research, Dr.\ Homi Bhabha Road, Pashan, Pune 411021, India.} 
\email{raghuram@iiserpune.ac.in}
\date{\today}

\thanks{This work is partially supported by an Alexander von Humboldt Research Fellowship.}


\maketitle


\def\g{\mathfrak{g}}
\def\k{\mathfrak{k}}
\def\z{\mathfrak{z}}
\def\h{{\mathfrak h}}
\def\gl{\mathfrak{gl}}

\def\q{{\tt q}}
\def\b{{\tt b}}
\def\t{{\tt t}}

\def\Ext{{\rm Ext}}
\def\Hom{{\rm Hom}}
\def\Ind{{\rm Ind}}

\def\GL{{\rm GL}}
\def\SL{{\rm SL}}
\def\SO{{\rm SO}}
\def\O{{\rm O}}

\def\R{\mathbb{R}}
\def\C{\mathbb{C}}
\def\Z{\mathbb{Z}}
\def\Q{\mathbb{Q}}
\def\A{\mathbb{A}}

\def\E{\mathcal{E}}
\def\G{\mathcal{G}}
\def\F{\mathcal{F}}
\def\T{\mathcal{T}}

\def\w{\wedge}

\def\Cat{\mathcal{C}}
\def\HC{{\rm HC}}
\def\HCat{\Cat^\HC}
\def\proj{{\rm proj}}
\def\sym{{\rm Sym}}

\def\autc{{\rm Aut}(\C)}

\def\to{\rightarrow}
\def\To{\longrightarrow}

\def\1{1\!\!1}
\def\dim{{\rm dim}}

\def\th{^{\rm th}}
\def\isom{\approx}

\def\CE{\mathcal{C}\mathcal{E}}


\numberwithin{equation}{section}
\newtheorem{thm}[equation]{Theorem}
\newtheorem{cor}[equation]{Corollary}
\newtheorem{lemma}[equation]{Lemma}
\newtheorem{prop}[equation]{Proposition}
\newtheorem{con}[equation]{Conjecture}
\newtheorem{ass}[equation]{Assumption}
\newtheorem{defn}[equation]{Definition}
\newtheorem{rem}[equation]{Remark}
\newtheorem{exer}[equation]{Exercise}
\newtheorem{exam}[equation]{Example}

\section{Introduction}
\label{sec:intro}

The Ramanujan $\Delta$-function is defined as 
$$
\Delta(z) \ = \  
q \prod_{n=1}^\infty (1-q^n)^{24} \ = \ 
\sum_{n=1}^\infty \tau(n) q^n, \quad q = e^{2 \pi i z}, \quad \Im(z) > 0.
$$
It is, up to scalar multiples, the unique weight $12$ holomorphic cusp form for $\SL_2(\Z).$ The Fourier coefficients 
$\tau(n)$ are multiplicative and hence the $L$-function of $\Delta$, defined initially as a Dirichlet series, admits an Euler product: 
$$
L(s, \Delta) \ := \  \sum_{n=1}^\infty \frac{\tau(n)}{n^s} \ = \ \prod_p (1 - \tau(p)p^{-s} + p^{11-2s})^{-1}, \quad \Re(s) \gg 0.
$$
The reciprocal of the Euler factor at $p$ may be factored as:
$$
1 - \tau(p)p^{-s} + p^{11-2s} \ = \ (1 - \alpha_p p^{-s})(1 - \beta_p p^{-s}). 
$$
For any $n \geq 1$, the $n$-th symmetric power $L$-function attached to $\Delta$ is defined as: 
$$
L(s, \sym^n, \Delta) \ = \ \prod_p \prod_{j=0}^n \left(1 - \alpha_p^j \beta_p^{n-j} p^{-s}\right)^{-1}, \quad \Re(s) \gg 0.
$$
The Langlands program predicts that this function admits an analytic continuation to an entire function, and satisfies a functional equation of the expected kind. This is known for $n \leq 4$ by the works of Hecke, Shimura, Gelbart and Jacquet, and Kim and Shahidi; furthermore, partial results are known for $5 \leq n \leq 9.$ 
See \cite{raghuram-shahidi} and the references therein. 

Consider now the symmetric cube $L$-function, i.e., take $n=3.$ 
We are interested in algebraicity results for the special values of 
$L(s, \sym^3, \Delta)$. The critical points for this $L$-function are $\{12, 13, 14, \dots, 22\}$; the center of symmetry for the functional equation here is $s = 17.$ There are several approaches to study these special values. 
Our general approach has been to use known cases of Langlands' principle of functoriality and the theory of cohomology of arithmetic groups to give a cohomological interpretation to certain analytic theories of $L$-functions. Let's adumbrate 
this theme: Let $\pi = \pi(\Delta)$ denote the cuspidal automorphic representation of $\GL_2(\A)$, where 
$\A$ is the ad\`ele ring of $\Q$. Appealing to the Langlands program, for any $n \geq 1$, let $\sym^n(\pi)$ denote the 
$n$-th symmetric transfer which is expected to be a cuspidal automorphic representation of $\GL_{n+1}(\A).$ 
The $n$-th symmetric transfer has been proven for $n \leq 4.$ The symmetric cube $L$-function of $\Delta$ appears as a factor in the Rankin-Selberg $L$-function for $\GL_3 \times \GL_2$, where on $\GL_3$ we take $\sym^2(\pi)$ and on $\GL_2$ we take $\pi$ itself; this approach was pursued in my paper \cite{raghuram-imrn} to get algebraicity results for 
$L(s, \sym^3(\pi))$, 
which is, up to an explicit shift in the $s$-variable, the same as $L(s, \sym^3, \Delta).$ Yet another way to get the special values is 
to observe that $\sym^3(\pi)$ has a so-called Shalika model and for such a representation an integral representation of the standard $L$-function $L(s, \sym^3(\pi))$ can be interpreted in cohomology (see my paper with 
Grobner \cite{grobner-raghuram-2}) giving special value results. Typically, one begins by identifying invariants or periods attached to representations such as $\sym^3(\pi),$ and then one proves via a cohomological interpretation of $L$-values, 
that $L(s, \sym^3(\pi))$ at a critical point $s = m$ is, up to rational numbers and up to archimedean factors (expected to be powers of $(2\pi i)$), given by one of these periods.  This begets the problem of comparing the various periods attached to a given representation. {\it The purpose of this article to provide comparison results between three different families of periods attached to a given representation.}

\medskip

We note that $ \sym^3(\pi)$ is a cuspidal automorphic representation of $\GL_4(\A)$ of cohomological type, i.e., 
contributes to the cohomology of a locally symmetric space with coefficients in a local system 
coming from an algebraic irreducible representation of the algebraic group $\GL_4/\Q$; for more details 
see \cite[Thm.\,5.5]{raghuram-shahidi}.
More generally, consider a cohomological cuspidal automorphic representation $\Pi$ of ${\rm GL}_{2n}(\A_F),$
where $F$ is a totally real number field.  We consider the following three kinds of periods attached to $\Pi$:

\smallskip

\begin{enumerate}
\item {\it Whittaker-Betti periods,} $p^\epsilon(\Pi_f)$, defined by comparing rational structures on two different models for the finite part $\Pi_f$ of the representation: one is a Whittaker model $W(\Pi_f)$ and the other is a cohomological model 
$H^\b(  \Pi_f  \times \epsilon);$ here $\epsilon$ is a $d$-tuple of signs where $d = [F:\Q]$ is the degree of $F$ over $\Q$. The cohomology degree $\b = dn^2$ is the lowest degree in which $\GL_{2n}/F$ has nontrivial cuspidal cohomology. 
Such periods, which can be defined in the context of cusp forms on $\GL_m$ over a number field, appear in the special values of Rankin-Selberg $L$-functions for $\GL_m \times \GL_{m-1}$. See my paper with Shahidi
\cite{raghuram-shahidi-imrn} together with my papers \cite{raghuram-imrn} and \cite{raghuram-preprint}. 

\smallskip

\item {\it Shalika-Betti periods,} $\omega^\epsilon(\Pi_f)$, defined by comparing rational structures on 
a Shalika model $S(\Pi_f)$, provided $\Pi$ does indeed admit a Shalika model, 
and a cohomological model $H^{\tt t}(\Pi_f \times \epsilon)$. 
The cohomology degree $\t = d(n^2+n-1)$ is the top-most degree in which $\GL_{2n}/F$ has nontrivial cuspidal cohomology. Such periods appear in algebraicity results for the critical values of the standard $L$-function for 
$\GL_{2n}/F$. See my paper with Grobner \cite{grobner-raghuram-2}. Let us note that ${\rm Sym}^3(\pi(\Delta))$, which is a representation of $\GL_4/\Q$, does indeed admit a Shalika model; see \cite[\S\,8.1]{grobner-raghuram-2}.

\smallskip

\item {\it Relative periods,} $\Omega_{\q}^{\epsilon}(\Pi_f)$ defined by comparing rational structures on two cohomological models $H^{\tt q}(\Pi_f \times \epsilon)$ and 
$H^{\tt q}(\Pi_f \times \epsilon')$ attached to $\Pi_f$; 
here $\q = \b$ the bottom degree or $\q = \t$ the top degree for nonvanishing cuspidal cohomology, and 
if $\epsilon = (\epsilon_v)_{v \in \Sigma_F}$ is a $d$-tuple of signs indexed by the set $\Sigma_F$ 
of all real embeddings of $F$, 
then $\epsilon' := -\epsilon = (-\epsilon_v)_{v \in \Sigma_F}$ is the opposite sign. 
Such periods, especially for $\q = \b$, appear in ratios of successive critical values for Rankin-Selberg $L$-functions for 
$\GL_{2n} \times \GL_{2m+1}.$ See my announcement with Harder \cite{harder-raghuram} and the forthcoming 
\cite{harder-raghuram-2}.
\end{enumerate}

\smallskip

In the main theorem of this article (see Thm.\,\ref{thm:main}) we prove two reciprocity laws. The first one relating 
Whittaker-Betti periods with relative periods in degree $\b$ says that the complex number 
$$
\frac{p^\epsilon(\Pi_f)}{p^{-\epsilon}(\Pi_f)} \, \frac{1}{\Omega^\epsilon_\b(\Pi_f)}
$$
is equivariant under the action of $\autc$. Since the representation $\Pi_f$ is defined over a number field, one may also say that this complex number is algebraic and is equivariant under  the  absolute Galois group over $\Q.$ This implies  that 
$p^\epsilon(\Pi_f)/p^{-\epsilon}(\Pi_f) 
 \ \approx \ \Omega^\epsilon_\b(\Pi_f),$
where, by $\approx$, we mean equality up an element of a number field $\Q(\Pi,\epsilon)$ 
which is the compositum of the rationality fields $\Q(\Pi)$ and $\Q(\epsilon)$ of 
$\Pi$ and $\epsilon$, respectively. 
The second reciprocity law looks similar and relates the Shalika-Betti periods with relative periods in degree $\t.$
We use these reciprocity laws to show that the relative periods are essentially invariant under twisting by characters; 
see Cor.\,\ref{cor} for a precise statement.

\smallskip

In \S\,\ref{sec:prelims} we briefly recall some basics about cohomological representations and the definitions of these 
families of periods. The main results are stated in \S\,\ref{sec:theorem} and their proofs are given in \S\,\ref{sec:proof}. The proof of Thm.\,\ref{thm:main} is modeled along the lines of the proof of the period relations as in my paper with Shahidi \cite{raghuram-shahidi-imrn}.  
One considers a certain diagram of isomorphisms of $\GL_{2n}(\A_{F,f})$-modules and the requisite reciprocity law captures the failure of commutativity of this diagram up to a corresponding diagram of rational structures; see the diagrams in (\ref{eqn:diagram-1}) and (\ref{eqn:diagram-2}).

\section{Notation and preliminaries}
\label{sec:prelims}

I will be very brief here, basically to set up the notations and define the periods. Any serious reader of this article is recommended to have by his/her side the following:  my paper with Freydoon~Shahidi \cite{raghuram-shahidi-imrn} for  the Whittaker-Betti periods; my paper with Harald~Grobner \cite{grobner-raghuram-2} together with my paper with Wee~Teck~Gan \cite{gan-raghuram} for the Shalika-Betti periods; and my paper with G\"unter~Harder \cite{harder-raghuram} for the relative periods.

\medskip
\subsection{Cohomological representations}
\label{sec:cohomology}
(Main Reference: \cite[\S\,3.3]{raghuram-shahidi}.) 
Let $F$ be a totally real field of degree $d = [F:\Q].$ We identify the set $S_\infty$ of infinite places of $F$ 
with the set $\Sigma_F = {\rm Hom}(F,\C) = 
{\rm Hom}(F, \R)$ of all real embeddings of $F$. 
Let $T$ be the diagonal torus in the standard Borel subgroup $B$ of upper triangular matrices in $\GL_{2n}/F.$ 
Let $\mu = (\mu_v)_{v \in \Sigma_F},$ a character of $T$, be a dominant integral weight. Let $(\rho_\mu, M_\mu)$ be the algebraic finite-dimensional irreducible representation of $\GL_{2n}(F \otimes \C) = \prod_{v \in S_\infty}\GL_{2n}(\C)$ with highest weight $\mu.$ 

\smallskip

A cuspidal automorphic representation $\Pi$ is said to be cohomological if there is a dominant integral weight 
$\mu$ such that the relative Lie algebra cohomology group
$$
H^\bullet(\g_\infty, K_\infty^\circ; \Pi_\infty \otimes M_\mu) \ \neq \ 0, 
$$
where, $\g_\infty$ is the complexified Lie algebra of $G_\infty = \GL_{2n}(F \otimes \R) = 
\prod_{v \in S_\infty} \GL_{2n}(\R);$ $K_\infty = \prod_{v \in S_\infty}({\rm O}(2n)Z_{2n}(\R))$ is the maximal compact subgroup times the center of $G_\infty$ and $K_\infty^\circ$ is the connected component of the identity in $K_\infty;$ 
hence $K_\infty^\circ = \prod_{v \in S_\infty}{\rm SO}(2n)Z_{2n}(\R)^\circ.$ For a cuspidal $\Pi$ as above, we will write 
$\Pi \in {\rm Cusp}(\GL_{2n}/F, \mu)$ to denote that it is cohomological with respect to $\mu$-coefficients. 

Given $\Pi \in {\rm Cusp}(\GL_{2n}/F, \mu)$,  it is known that 
$$
H^\bullet(\g_\infty, K_\infty^\circ; \Pi_\infty \otimes M_\mu) \neq 0 
\  \Longleftrightarrow \ 
dn^2 \leq \bullet \leq d(n^2+n-1).
$$
Let $\b = dn^2$ denote the bottom-degree, and $\t = d(n^2+n-1)$ the top-degree for non-vanishing of the above cohomology group. Let $\q = \b$ or $\q = \t$ be either of these two extreme degrees. 
The group of connected components
$\pi_0(K_\infty)$ acts on these  relative Lie algebra cohomology groups. 
As a $\pi_0(K_\infty) \cong \prod_{v \in S_\infty} {\rm O}(2n)/{\rm SO}(2n)$-module we have:
$$
H^\q(\g_\infty, K_\infty^\circ; \Pi_\infty \otimes M_\mu) \ = \ 
\bigoplus\limits_{\epsilon \in \widehat{\pi_0(K_\infty)}} \epsilon.
$$
Let $V_\Pi$ denote the unique subspace of the space of cusp forms on $\GL_{2n}(\A_F)$ realizing the 
representation $\Pi.$ Then in degree $\q = \b$ or $\q = \t$ we have 
$$
H^\q(\g_\infty, K_\infty^\circ; V_\Pi \otimes M_\mu) \ = \ 
\bigoplus\limits_{\epsilon \in \widehat{\pi_0(K_\infty)}} \epsilon \otimes \Pi_f
$$
as a $\pi_0(K_\infty) \times \GL_{2n}(\A_{F,f})$-module; where $\A_{F,f}$ is the ring of finite ad\`eles of $F$. 
Therefore for each $\epsilon  \in \widehat{\pi_0(K_\infty)}$ and an extreme degree $\q$, the summand 
$ \epsilon \otimes \Pi_f$ in the right hand side of the above decomposition 
gives  a {\it cohomological model} 
for the finite part of the representation which we denote:
$$
H^\q(\Pi_f \times \epsilon).
$$

There is an action of $\autc$ on all these cohomological realizations. First of all, $\autc$ acts on $\Sigma_F$ by  composition, or one can say $\autc$ acts on $S_\infty$ by permutations induced by composition. 
Next, given a cuspidal representation $\Pi$, define an abstract representation 
${}^\sigma\Pi = {}^\sigma\Pi_\infty \otimes {}^\sigma\Pi_f$ as in 
\cite[\S\,2.5]{grobner-raghuram-2}. It is deep theorem of Clozel \cite[Thm.\,3.13]{clozel} 
that if $\Pi \in {\rm Cusp}(\GL_{2n}/F, \mu)$ then  
${}^\sigma\Pi \in {\rm Cusp}(\GL_{2n}/F, {}^\sigma\mu),$ and we have a $\sigma$-semi-linear, $\GL_{2n}(\A_{F,f})$-equivariant isomorphism:
$$
\tilde\sigma : H^\q(\Pi_f \times \epsilon) \longrightarrow H^\q({}^\sigma\Pi_f \times {}^\sigma\epsilon), \quad \q \in \{\b, \t \}.
$$

Define the rationality field $\Q(\Pi_f)$ as the subfield of $\C$ fixed by 
$\{\sigma \in \autc : {}^\sigma\Pi_f \simeq \Pi_f\}.$ Similarly, define $\Q(\mu)$ as  the subfield of $\C$ fixed by 
$\{\sigma \in \autc : {}^\sigma\mu= \mu\}.$ By strong multiplicity one for $\GL_n$, we have 
$\Q(\mu) \subset \Q(\Pi_f).$ Now rename $\Q(\Pi_f)$ as $\Q(\Pi).$   
For any $\epsilon$ as above, define $\Q(\epsilon)$ in a similar way; and $\Q(\Pi, \epsilon)$ is the compositum of 
$\Q(\Pi)$ and $\Q(\epsilon).$ We also know from Clozel, {\it loc.cit.,} that the cohomological models 
$H^\q(\Pi_f \times \epsilon)$ have a $\Q(\Pi, \epsilon)$ structure.

\subsection{Whittaker-Betti periods}
(Main Reference:  \cite{raghuram-shahidi-imrn}.) 
Fix a non-trivial additive character $\psi : F\backslash \A_F \to \C^\times$ as, for example, in Tate's thesis. 
A cuspidal automorphic representation $(\Pi,V_\Pi)$ of $\GL_{2n}(\A_F)$ has a global $\psi$-Whittaker model, and hence a local $\psi_v$-Whittaker model for all places $v$ of $F.$  This gives a 
Whittaker model $W(\Pi_f) := \otimes_{v \notin S_\infty} W(\Pi_v)$ for  $\Pi_f$. 
For $\sigma \in \autc,$ as described in \cite[\S\,3.2]{raghuram-shahidi-imrn}, we get a 
$\sigma$-semi-linear, $\GL_{2n}(\A_{F,f})$-equivariant isomorphism:
$$
\tilde\sigma : W(\Pi_f) \longrightarrow W({}^\sigma\Pi_f).
$$
In particular, this says that $W(\Pi_f)$ has a $\Q(\Pi_f)$-structure which we denote $W(\Pi_f)_0$; one can see that a rational structure is generated by a new-vector; see \cite[Lem.\,3.2]{raghuram-shahidi-imrn}. 

\smallskip

There is a comparison isomorphism between the Whittaker model $W(\Pi_f)$ and a cohomological model 
$H^\b(\Pi_f \times \epsilon),$ for any $\epsilon  \in \widehat{\pi_0(K_\infty)}.$ 
This isomorphism depends on a choice of a generator $[\Pi_\infty]^\epsilon$ of the one-dimensional subspace of 
$H^\b(\g_\infty, K_\infty^\circ; \Pi_\infty \otimes M_\mu)$  realizing $\epsilon.$ These generators are chosen in a compatible fashion, i.e., $[\Pi_\infty]^\epsilon$ and $[{}^\sigma\Pi_\infty]^{{}^\sigma\epsilon}$ are related by the permutation action of $\sigma$ on $S_\infty.$ 
We take these generators to be  fixed for the rest of the paper. The comparison isomorphism is denoted by: 
$$
\F^\epsilon(\Pi_f) : W(\Pi_f) \longrightarrow H^\b(\Pi_f \times \epsilon). 
$$
The rational structure on $W(\Pi_f)$ has no relation to the rational structure on $H^\b(\Pi_f \times \epsilon).$ However, both sides are representation spaces for an irreducible action of $\GL_{2n}(\A_{F,f})$. There is a uniqueness principle at play in such a situation (see Waldspurger \cite[Lem.\,I.1]{waldy}) which defines certain periods: there exists 
$p^\epsilon(\Pi_f) \in \C^\times$ such that 
\begin{equation}
\label{eqn:wb}
\tilde\sigma \circ \left( \frac{1}{p^\epsilon(\Pi_f)}  \F^\epsilon(\Pi_f) \right) \ = \ 
\frac{1}{ p^{{}^\sigma\!\epsilon}({}^\sigma\Pi_f)}  \F^{{}^\sigma\!\epsilon}({}^\sigma\Pi_f) \circ \tilde\sigma;  
\end{equation}
i.e., the normalized map $\frac{1}{p^\epsilon(\Pi_f)}  \F^\epsilon(\Pi_f)$ preserves rational structures. 
See \cite[Def./Prop.\,3.3]{raghuram-shahidi-imrn}.

\subsection{Shalika-Betti periods}
(Main reference: \cite{grobner-raghuram-2}.)
Suppose now that the representation $\Pi \in {\rm Cusp}(\GL_{2n}/F, \mu)$ admits a Shalika model; 
see \cite[\S\,3]{grobner-raghuram-2}. Any discussion of a Shaika model for $\Pi$ entails a choice of an additive character $\psi$ as with Whittaker models, and a choice of a Hecke character $\eta : F^\times\backslash \A_F^\times \to \C^\times$ such that $\eta^n = \omega_\Pi$ the central character of $\Pi.$ We may take $\eta$ to be an algebraic Hecke character and then we fix our choice of $\eta$; we also have its rationality field $\Q(\eta).$
To emphasize the dependence on these characters, we say that $\Pi$ has an $(\eta, \psi)$-Shalika model. 
(See \cite[Rem.\,3.8.2]{grobner-raghuram-2} on multiplicity one for Shalika models.) 
For $\sigma \in \autc$,  the representation 
${}^\sigma\Pi$ has a $({}^\sigma\eta,\psi)$-Shalika model by the arithmeticity property for Shalika periods; 
see \cite{gan-raghuram}. As described in \cite[\S\,3.7]{grobner-raghuram-2} we have a 
$\sigma$-semi-linear, $\GL_{2n}(\A_{F,f})$-equivariant isomorphism between the $(\eta,\psi)$-Shalika model of 
$\Pi_f,$ denoted for brevity as $S(\Pi_f),$ and the $({}^\sigma\eta, \psi)$-Shalika model 
$S({}^\sigma\Pi_f)$ of ${}^\sigma\Pi_f$: 
$$
\tilde\sigma : S(\Pi_f) \longrightarrow S({}^\sigma\Pi_f).
$$
In particular, this says that $S(\Pi_f)$ has a $\Q(\Pi_f, \eta)$-structure which we denote $S(\Pi_f)_0;$ furthermore, one can see that $S(\Pi_f)_0$ is generated by a new-vector; see \cite[Lem.\,3.8.1]{grobner-raghuram-2}. 
For $\epsilon \in \widehat{\pi_0(K_\infty)}$ we fix a generator (which are all compatible as before) for the 
the one-dimensional subspace of 
$H^\t(\g_\infty, K_\infty^\circ; \Pi_\infty \otimes M_\mu)$  realizing $\epsilon.$ 
We have a comparison isomorphism (see \cite[\S\,4.2]{grobner-raghuram-2}): 
$$
\Theta^\epsilon(\Pi_f) :  S(\Pi_f) \longrightarrow H^\t(\Pi_f \times \epsilon). 
$$
The rational structure on $S(\Pi_f)$ has no relation to the rational structure on $H^\t(\Pi_f \times \epsilon),$ and this situation defines $\omega^\epsilon(\Pi_f) \in \C^\times$ such that 
\begin{equation}
\label{eqn:sb}
\tilde\sigma \circ \left( \frac{1}{\omega^\epsilon(\Pi_f)}  \Theta^\epsilon(\Pi_f) \right) \ = \ 
\frac{1}{ \omega^{{}^\sigma\!\epsilon}({}^\sigma\Pi_f)}  \Theta^{{}^\sigma\!\epsilon}({}^\sigma\Pi_f) \circ \tilde\sigma;  
\end{equation}
i.e., the normalized map $\frac{1}{\omega^\epsilon(\Pi_f)}  \Theta^\epsilon(\Pi_f)$ preserves rational structres. 
See \cite[Def./Prop.\,4.2.1]{grobner-raghuram-2}.

\subsection{Relative periods}
\label{sec:rel-periods}
(Main Reference: \cite[\S\,2]{harder-raghuram}.) 
These periods are intrinsic to the theory of cohomology of arithmetic groups and their definition is independent of any transcendentally defined models like the Whittaker or the Shalika models considered above. 
Given $\epsilon$ and its opposite sign $\epsilon'$ we have an isomorphism of $\GL_{2n}(\A_{F,f})$-modules:  
\begin{equation}\label{eqn:relative-iso}
\T_\q^\epsilon(\Pi_f) : H^\q(\Pi_f \times \epsilon) \longrightarrow H^\q(\Pi_f \times \epsilon'), \quad \q \in \{\b,\t\}. 
\end{equation}
(See \cite[(2)]{harder-raghuram}.) Later, somewhat artificially, we will arrange for this isomorphism to be related to an appropriate (Whittaker or Shalika) model; such an arrangement is what will ultimately make the proof work; 
see (\ref{eqn:triangle-w}) and (\ref{eqn:triangle-s}) below. 
There is no reason for this isomorphism to preserve rational structures on either side, and this defines the relative 
period $\Omega_\q^\epsilon(\Pi_f) \in \C^\times$ as that complex number needed to modify the isomorphism to make it $\autc$-equivariant: 
\begin{equation}
\label{eqn:relative}
\tilde\sigma \circ \left( \Omega_\q^\epsilon(\Pi_f) \,  \T_\q^\epsilon(\Pi_f) \right) \ = \ 
\Omega_\q^{{}^\sigma\!\epsilon}({}^\sigma\Pi_f) \, \T_\q^{{}^\sigma\!\epsilon}({}^\sigma\Pi_f) \circ \tilde\sigma. 
\end{equation}

\medskip
\section{Statement of the main results on period relations}
\label{sec:theorem}

\smallskip

\begin{thm}
\label{thm:main}
Let $\Pi$ be a cohomological cuspidal automorphic representation of $\GL_{2n}$ over a totally real number field $F$. 
\begin{enumerate}
\item The Whittaker-Betti periods $p^\epsilon(\Pi_f)$ 
are related to the relative periods $\Omega^\epsilon_\b(\Pi_f)$ in degree $\b$ as follows: 
$$
\sigma\left(
\frac{p^\epsilon(\Pi_f)}{p^{-\epsilon}(\Pi_f)} \,  \frac{1}{\Omega^\epsilon_\b(\Pi_f)} \right) 
\ = \ 
\frac{p^{{}^\sigma\epsilon}({}^\sigma\Pi_f)}{p^{-{}^\sigma\epsilon}({}^\sigma\Pi_f)}  \, 
\frac{1}{\Omega^{{}^\sigma\epsilon}_\b({}^\sigma\Pi_f)},  \quad \forall \sigma \in \autc.
$$
In particular, we have 
$$
\frac{p^\epsilon(\Pi_f)}{p^{-\epsilon}(\Pi_f)} 
 \ \approx \ \Omega^\epsilon_\b(\Pi_f), 
$$
where, by $\approx$, we mean equality up to an element of the number field $\Q(\Pi, \epsilon)$. 

\medskip

\item Suppose that $\Pi$ admits an $(\eta, \psi)$-Shalika model. Then the Shalika-Betti periods $\omega^\epsilon(\Pi_f)$
are related to the relative periods $\Omega^\epsilon_\t(\Pi_f)$ in degree $\t$ as follows: 
$$
\sigma\left(
\frac{\omega^\epsilon(\Pi_f)}{\omega^{-\epsilon}(\Pi_f)} \, 
\frac{1}{\Omega^\epsilon_\t(\Pi_f)} \right) 
\ = \ 
\frac{\omega^{{}^\sigma\epsilon}({}^\sigma\Pi_f)}{\omega^{-{}^\sigma\epsilon}({}^\sigma\Pi_f)} \, 
\frac{1}{\Omega^{{}^\sigma\epsilon}_\t({}^\sigma\Pi_f)}, \quad \forall \sigma \in \autc.
$$
In particular, we have 
$$
\frac{\omega^\epsilon(\Pi_f)}{\omega^{-\epsilon}(\Pi_f)} 
 \ \approx \ \Omega^\epsilon_\t(\Pi_f), 
$$
where, by $\approx$, we mean equality up to an element of the number field $\Q(\Pi, \epsilon, \eta)$. 
\end{enumerate}
\end{thm}

\begin{cor}[Behaviour of relative periods under twisting by characters]
\label{cor}
Let $\Pi$ be a cohomological cuspidal automorphic representation of $\GL_{2n}$ over a totally real number field $F$ and let $\chi$ be an algebraic Hecke character of $F$. 
Then: 
\begin{enumerate}
\item For $\sigma \in \autc$ we have
$$
\sigma\left(\frac{\Omega_\b^\epsilon(\Pi_f \otimes \chi_f)}
{\Omega_\b^{\epsilon\epsilon_\chi}(\Pi_f)}\right) \ = \ 
\frac{\Omega_\b^{{}^\sigma\epsilon}({}^\sigma\Pi_f \otimes {}^\sigma\chi_f)}
{\Omega_\b^{{}^\sigma\epsilon \, {}^\sigma\epsilon_\chi}({}^\sigma\Pi_f)}. 
$$
In particular, we have 
$$
\Omega_\b^\epsilon(\Pi_f \otimes \chi_f) \ \approx \ 
\Omega_\b^{\epsilon\epsilon_\chi}(\Pi_f), 
$$
where, by $\approx$, we mean equality up to an element of $\Q(\Pi, \epsilon, \chi)$.

\smallskip

\item For $\sigma \in \autc$ we have 
$$
\sigma\left(\frac{\Omega_\t^\epsilon(\Pi_f \otimes \chi_f)}
{\Omega_\t^{\epsilon\epsilon_\chi}(\Pi_f)}\right) \ = \ 
\frac{\Omega_\t^{{}^\sigma\epsilon}({}^\sigma\Pi_f \otimes {}^\sigma\chi_f)}
{\Omega_\t^{{}^\sigma\epsilon \, {}^\sigma\epsilon_\chi}({}^\sigma\Pi_f)}. 
$$
(We do not need to assume that $\Pi$ has a Shalika model.) In particular, we have 
$$
\Omega_\t^\epsilon(\Pi_f \otimes \chi_f) \ \approx \ 
\Omega_\t^{\epsilon\epsilon_\chi}(\Pi_f), 
$$
where, by $\approx$, we mean equality up to an element of $\Q(\Pi, \epsilon, \chi)$.
\end{enumerate}
\end{cor}

\bigskip
\section{Proofs}
\label{sec:proof}

Let's begin with the proof of (1) of Theorem~\ref{thm:main} which is a comparison between the 
Whittaker-Betti periods $p^\epsilon(\Pi_f)$ and the relative periods $\Omega^\epsilon_\b(\Pi_f)$ in degree $\b.$ 
The proof is roughly along the lines of the proof of Theorem 4.1 in my paper with 
Shahidi \cite{raghuram-shahidi-imrn}. Consider the following diagram of isomorphisms of $\GL_{2n}(\A_{F,f})$-modules. 
{\it Warning:} This diagram is not commutative! Indeed the period relation captures the failure of commutativity up to a corresponding diagram of rational structures. 

{\SMALL
\begin{equation}\label{eqn:diagram-1}
\xymatrix{
H^\b(\Pi_f \times \epsilon) \ar[rr]^{\T_\b^\epsilon(\Pi_f)} \ar[rrrd]^{\tilde\sigma} & 
& H^\b(\Pi_f \times \epsilon') \ar[rrrd]^{\tilde\sigma} & & & \\
 & & & H^\b({}^\sigma\Pi_f \times {}^\sigma\!\epsilon) \ar[rr]_{\T_\b^{{}^\sigma\!\epsilon}({}^\sigma\Pi_f)} & 
 & H^\b({}^\sigma\Pi_f \times {}^\sigma\!\epsilon')  \\
 & W(\Pi_f) \ar[luu]^{\F^\epsilon(\Pi_f)} \ar[ruu]_{\F^{\epsilon'}\!(\Pi_f)} \ar[rrrd]^{\tilde\sigma} & & & & \\
 & & & & W({}^\sigma\Pi_f) \ar[luu]^{\F^{{}^\sigma\!\epsilon}({}^\sigma\Pi_f)} 
 \ar[ruu]_{\F^{{}^\sigma\!\epsilon'}\!({}^\sigma\Pi_f)} & 
}
\end{equation}}

\noindent
In this diagram, consider the composition $\tilde\sigma \circ \T_\b^\epsilon(\Pi_f) \circ \F^\epsilon(\Pi_f)$. We are going to evaluate this in two different ways. 

On the one hand, using (\ref{eqn:relative}) for the relative periods, and then using 
(\ref{eqn:wb}) for the Whittaker-Betti periods, we have
\begin{equation}
\label{eqn:step1}
\begin{array}{lll}
\tilde\sigma \circ \T_\b^\epsilon(\Pi_f) \circ \F^\epsilon(\Pi_f) & = & 
\frac{\Omega_\b^{{}^\sigma\!\epsilon}({}^\sigma\Pi_f)}{\sigma(\Omega_\b^\epsilon(\Pi_f))} \, 
\T_\b^{{}^\sigma\!\epsilon}({}^\sigma \Pi_f) \circ \tilde\sigma \circ  \F^\epsilon(\Pi_f) \\ 
 & & \\ 
 & = & 
 \frac{\Omega_\b^{{}^\sigma\!\epsilon}({}^\sigma\Pi_f)}{\sigma(\Omega_\b^\epsilon(\Pi_f))} \, 
\frac{\sigma(p^\epsilon(\Pi_f))}{p^{{}^\sigma\!\epsilon}({}^\sigma\Pi_f)} \, 
\T_\b^{{}^\sigma\!\epsilon}({}^\sigma \Pi_f) \circ  \F{{}^\sigma\!\epsilon}({}^\sigma \Pi_f) \circ \tilde\sigma. 
\end{array}
\end{equation}

Now we come to an important point that the isomorphism $\T_\b^\epsilon(\Pi_f)$ in (\ref{eqn:relative-iso}), needed to define the relative periods in bottom degree, may be taken to be so that
\begin{equation}\label{eqn:triangle-w}
\T_\b^\epsilon(\Pi_f) \circ \F^\epsilon(\Pi_f) =  \F^{\epsilon'}(\Pi_f).
\end{equation}
This follows from \cite[(2)]{harder-raghuram}. (In the notations of \cite{harder-raghuram}, the Hecke summand $\pi_f$ appearing rationally in inner cohomology may be identified with a rational structure $W(\Pi_f)_0$ on the Whittaker model; after making this identification over $\Q(\Pi_f)$, (\ref{eqn:triangle-w}) follows from  \cite[(2)]{harder-raghuram}.)
In other words, if the diagram~(\ref{eqn:diagram-1}) is thought of as a triangular prism, then the triangles at both the ends commute. 

So, on the other hand, using (\ref{eqn:triangle-w}) and (\ref{eqn:wb}), we have 
\begin{equation}
\label{eqn:step2}
\begin{array}{lll}
\tilde\sigma \circ \T_\b^\epsilon(\Pi_f) \circ \F^\epsilon(\Pi_f) & = & \tilde\sigma \circ \F^{\epsilon'}(\Pi_f) \\
& & \\ 
 & = & 
 \frac{\sigma(p^{\epsilon'}(\Pi_f))}{p^{{}^\sigma\!\epsilon'}({}^\sigma\Pi_f)} \, 
 \F^{{}^\sigma\!\epsilon'}({}^\sigma\Pi_f) \circ \tilde\sigma \\ 
 & & \\
 & = & 
  \frac{\sigma(p^{\epsilon'}(\Pi_f))}{p^{{}^\sigma\!\epsilon'}({}^\sigma\Pi_f)} \, 
  \T_\b^{{}^\sigma\!\epsilon}({}^\sigma \Pi_f) \circ  \F{{}^\sigma\!\epsilon}({}^\sigma \Pi_f) \circ \tilde\sigma. 
\end{array}
\end{equation}
Comparing the right hand sides of (\ref{eqn:step1}) and (\ref{eqn:step2}) we get the following relation: 
$$
 \frac{\Omega_\b^{{}^\sigma\!\epsilon}({}^\sigma\Pi_f)}{\sigma(\Omega_\b^\epsilon(\Pi_f))} \, 
\frac{\sigma(p^\epsilon(\Pi_f))}{p^{{}^\sigma\!\epsilon}({}^\sigma\Pi_f)}  \ = \ 
 \frac{\sigma(p^{\epsilon'}(\Pi_f))}{p^{{}^\sigma\!\epsilon'}({}^\sigma\Pi_f)}, 
$$
rewriting which proves (1) of Theorem~\ref{thm:main}. 

\medskip

The proof of (2) of Theorem~\ref{thm:main}  is almost identical. Here we consider the  diagram which is a triangular prism of $\GL_{2n}(\A_{F,f})$-module isomorphisms: 
{\SMALL
\begin{equation}\label{eqn:diagram-2}
\xymatrix{
H^\t(\Pi_f \times \epsilon) \ar[rr]^{\T_\t^\epsilon(\Pi_f)} \ar[rrrd]^{\tilde\sigma} & 
& H^\t(\Pi_f \times \epsilon') \ar[rrrd]^{\tilde\sigma} & & & \\
 & & & H^\t({}^\sigma\Pi_f \times {}^\sigma\!\epsilon) \ar[rr]_{\T_\t^{{}^\sigma\!\epsilon}({}^\sigma\Pi_f)} & 
 & H^\t({}^\sigma\Pi_f \times {}^\sigma\!\epsilon')  \\
 & S(\Pi_f) \ar[luu]^{\Theta^\epsilon(\Pi_f)} \ar[ruu]_{\Theta^{\epsilon'}\!(\Pi_f)} \ar[rrrd]^{\tilde\sigma} & & & & \\
 & & & & S({}^\sigma\Pi_f) \ar[luu]^{\Theta^{{}^\sigma\!\epsilon}({}^\sigma\Pi_f)} 
 \ar[ruu]_{\Theta^{{}^\sigma\!\epsilon'}\!({}^\sigma\Pi_f)} & 
}
\end{equation}}

\noindent
Now the composition $\tilde\sigma \circ \T_\t^\epsilon(\Pi_f) \circ \Theta^\epsilon(\Pi_f)$ is written in two different ways, 
again while using the observation that the triangles at either end are commutative. This needs a little explanation: 
In the isomorphism (\ref{eqn:relative-iso}) for top-degree relative periods, we can arrange it so that 
\begin{equation}\label{eqn:triangle-s}
\T_\t^\epsilon(\Pi_f) \circ \Theta^\epsilon(\Pi_f) =  \Theta^{\epsilon'}(\Pi_f).
\end{equation}
This follows from \cite[(2)]{harder-raghuram}. (But this time, in the notations of \cite{harder-raghuram}, the Hecke summand 
$\pi_f$ appearing rationally in inner cohomology is identified with a rational structure $S(\Pi_f)_0$ 
for the Shalika model; after making this identification over $\Q(\Pi_f)$,  (\ref{eqn:triangle-s}) follows from  
\cite[(2)]{harder-raghuram}.) We omit the rest of the details since they are absolutely similar to the proof of (1).

\medskip
For the proof of (1) of Corollary~\ref{cor}, write the ratio 
$$
\frac{\Omega_\b^\epsilon(\Pi_f \otimes \chi_f)}{\Omega_\b^{\epsilon\epsilon_\chi}(\Pi_f)}
$$
as a product of ratios: 
{\Small
$$
\left(\frac{\Omega_\b^\epsilon(\Pi_f \otimes \chi_f)p^{-\epsilon}(\Pi_f \otimes\chi_f)}{p^{\epsilon}(\Pi_f \otimes\chi_f)} \right)
\left(\frac{p^{\epsilon}(\Pi_f \otimes\chi_f)}{p^{\epsilon\epsilon_\chi}(\Pi_f)\G(\chi_f)^{\frac{n(n-1)}{2}}}\right)
\left(\frac{p^{-\epsilon \epsilon_\chi}(\Pi_f) \G(\chi_f)^{\frac{n(n-1)}{2}} } 
{p^{-\epsilon}(\Pi_f \otimes \chi_f)} \right)
\left(\frac{p^{\epsilon \epsilon_\chi}(\Pi_f)}{p^{-\epsilon \epsilon_\chi}(\Pi_f)\Omega_\b^{\epsilon\epsilon_\chi}(\Pi_f)}
\right), 
$$}

\noindent
where $\G(\chi_f)$ is the Gau\ss~sum of $\chi_f$ as defined in \cite[\S\,2]{raghuram-imrn}. 
The first and last factors are $\autc$-equivariant by (1) of Theorem~\ref{thm:main}, and the middle two factors are 
$\autc$-equivariant by the main theorem of my paper with Shahidi \cite{raghuram-shahidi-imrn} on the behavior of Whittaker-Betti periods under twisting by algebraic Hecke characters, but as rewritten in \cite[(2.2)]{raghuram-imrn}

\medskip

For the proof of (2) of Corollary~\ref{cor}, we could rewrite the ratio 
$\Omega_\t^\epsilon(\Pi_f \otimes \chi_f)/\Omega_\t^{\epsilon\epsilon_\chi}(\Pi_f)$ as above in terms of the Shalika-Betti periods and then use  \cite[Thm.\,5.2.1]{grobner-raghuram-2} about the behavior of these latter periods under twisting. However, we can also do this without recourse to Shalika models. First of all we define top-degree Whittaker-Betti periods,
denoted $q^\epsilon(\Pi_f)$, by comparing $W(\Pi_f)$ with $H^\t(\Pi_f \times \epsilon)$; see the last paragraph of \cite{raghuram-shahidi-imrn}. Next, we prove a variation of (1) of Theorem~\ref{thm:main}, via an identical proof, to show that 
\begin{equation}\label{eqn:variation}
\frac{q^\epsilon(\Pi_f)}{q^{-\epsilon}(\Pi_f)}\frac{1}{\Omega_\t^\epsilon(\Pi_f)}
\end{equation}
is $\autc$-equivariant. Now (2) of Corollary~\ref{cor} follows, exactly as in the proof of (1) of this corollary, by using 
(\ref{eqn:variation}) and 
\cite[(4.6)]{raghuram-shahidi-imrn} on the behavior of these $q^\epsilon(\Pi_f)$ under twisting. 
This concludes the proofs. We end this article with some speculative remarks. 

\begin{rem}{\rm 
I believe there should be some interesting relation between the relative periods 
$\Omega_\b^\epsilon(\Pi_f)$ in bottom-degree and the relative periods $\Omega_\t^\epsilon(\Pi_f)$ in top-degree,  possibly replacing one of the $\Pi$ by its contragredient representation $\Pi^{\sf v}$, and that such a relation should follow from Poincar\'e duality. However, I have not been able to make this precise.
}\end{rem}

\begin{rem}{\rm 
\label{rem:referee-comment}
The referee asked an interesting question as to why in the definition of relative periods (see \ref{sec:rel-periods}) one needs to work with a comparison of cohomological models with respect to a sign $\epsilon$ and its opposite sign 
$\epsilon' = -\epsilon.$ One can seemingly define relative periods, say $\Omega^{\epsilon_1,\epsilon_2}(\Pi_f)$,  
by playing off two cohomological models for $\Pi_f$ using any pair of signs $\epsilon_1$ and $\epsilon_2$. For such generalized relative periods, one can prove analogous period-relations as in Theorem~\ref{thm:main}. However, to the best of my knowledge, it is only when  we work with opposite signs is there a relation to $L$-values, and this is one of the main themes in my work with Harder \cite{harder-raghuram} \cite{harder-raghuram-2}. See also the discussion in the introduction of  my paper with Bhagwat \cite{bhagwat-raghuram} where it becomes clear that if we look at ratios of successive critical values then one is forced to be looking a sign and its opposite sign. It could be interesting to look for some arithmetic interpretation of these generalized relative periods $\Omega^{\epsilon_1,\epsilon_2}(\Pi_f)$. 
}\end{rem}

\medskip{\small
\noindent{\it Akcnowledgements:} It is a pleasure to thank Harald Grobner who made some interesting comments on a first draft of this article. I thank the organizers of the Legacy of Ramanujan Conference for inviting me to give a talk there. Finally, I acknowledge the pleasant working environment of the Max Planck Institut f\"ur Mathematik, Bonn, where this paper was written-up.}

\end{document}